\documentclass[11pt]{amsart}
\usepackage{amssymb}

\theoremstyle{plain}
\newtheorem{theorem}{Theorem}

\newtheorem{lemma}{Lemma}
\newtheorem{proposition}{Proposition}

\theoremstyle{definition}

\theoremstyle{remark}

\numberwithin{equation}{section}

\begin{document}
\title[Small Prime Powers in the Fibonacci Sequence ]
       {\Large{Small Prime Powers in the Fibonacci Sequence }}
\author{J. Mc Laughlin}
\address{Mathematics Department\\
       University of Illinois \\
        Champaign - Urbana, Illinois 61820}
\email{jgmclaug@math.uiuc.edu}
\date{December,13,2000}
\begin{abstract} 
It is shown that there are no non-trivial fifth-, seventh-,
eleventh-, thirteenth- or seventeenth powers in the Fibonacci
sequence.

For eleventh, thirteenth- and seventeenth powers an alternative 
 (to the usual exhaustive check of products of powers of fundamental
units) method is used to 
overcome the problem of having a large number of independent units and
relatively high bounds on their exponents.

It is envisaged that the same method can be used to decide the
question of the existence of higher small prime powers in the
Fibonacci sequence and that the method can be applied to 
  other binary recurrence sequences. The alternative method
 mentioned may have
wider applications.
\end{abstract}

\maketitle

\section{Introduction }
The Fibonacci sequence \{$F_{n}$\}$_{n=0}^{\infty}$ is defined by
setting $F_{0} = 0, F_{1} = 1$ and, for $n \geq 2$, by setting 
$F_{n} = F_{n-1} + F_{n-2}$ .

	Cohn~\cite{C64} and Wylie~\cite{W64} proved independently, by
elementary means, that
the only squares in the Fibonacci sequence are $F_{0}=0,F_{1}= F_{2} =
1$ and $F_{12} = 144$.

	In~\cite{LF70}, London and Finkelstein used previous results on
solutions to two diophantine equations to show that the only cubes in
the Fibonacci sequence are $F_{0}=0$, $F_{1}= F_{2} =
1$ and $F_{6} = 8$. In \cite{LW81}, Lagarias and Weisser 
gave a complete determination of all Fibonacci numbers of the form
$2^{a}3^{b}x^{3}$.
 In~\cite{P83}, Peth\H{o} used linear forms in
logarithms together with a computer search using congruence
considerations to give an alternative proof of London and
Finkelstein's result.

 Peth\H{o} (\cite{P82}) and Shorey and Stewart (\cite{SS83}) 
proved independently  that there are only finitely
many perfect powers in any non-trivial binary recurrence 
sequence.
	In~\cite{P81} Peth\H{o} states that
if $F_{m} = x^{q}$, for some positive
integers $x$ and $q$, then $q<10^{98}$.
  In the same
paper he also states that he used the same method
that he used in~\cite{P83} to show that the only fifth powers in the
Fibonacci sequence are $F_{0}=0$ and $F_{1}= F_{2} =1$.

	In this paper the method outlined by Peth\H{o} in \cite{P83}
is used as a starting point and then linear forms in logarithms together 
with the LLL
algorithm are used to reprove the result for  fifth powers
and to prove that the only seventh-,
eleventh-, thirteenth- or seventeenth powers in the Fibonacci sequence
are $F_{0}=0$ and $F_{1}= F_{2} =1$. An alternative method 
(to the usual exhaustive check of products of powers of fundamental
units)  is used to complete the
search in the case of eleventh, thirteenth- and seventeenth powers.

It is envisaged that the same method can be used to decide the
question of the existence of higher small prime powers in the
Fibonacci sequence and that the method can be applied to search for
prime powers in other binary reurrence sequences. The alternative method 
 mentioned above may have wider applications.

Computations were performed using {\itshape Magma\/}, 
{\itshape Mathematica\/} and 
{\itshape Pari-gp\/} and were carried out on Sun ultra 5- and Sun ultra 10 
computers.

\section{Elementary Considerations}

It is easy to show that 
{\allowdisplaybreaks
\begin{align*}
&\phantom{as}\\
F_{m+1}^{2}-F_{m+1}F_{m}-F_{m}^2 &= (-1)^{m}.\\
&\phantom{as}
\end{align*}
}
Peth\H{o}, in~\cite{P83},  gives the following lemma:
\begin{lemma}
Let $q \geq 3$ be a positive integer. If $F_{m} = x^{q}$, for some
positive integer $x$, then $m=0,1,2,6$ or $\exists$ a prime $p$,
$p|m$, 
such that $F_{p} = x_{1}^{q}$, for some positive integer $x_{1}$.
\end{lemma}
Hence it can be assumed that $m$ is an odd prime. Suppose
$F_{m}=x^{n}$ and $F_{m+1} = y$, for some positive integers $x,y$ and 
some prime 
$n \geq 5$. Then $y^{2}-yx^{n}-x^{2n}+1=0$, and regarding this equation as a 
quadratic in $y$ and looking at its discriminant, it follows that 
$5x^{2n}-4=z^{2}$, for some positive integer $z$. Clearly $z=5v \pm
1$, for some positive integer $v$ and thus that
{\allowdisplaybreaks
\begin{align}\label{eq:2}
&\phantom{as}\notag\\
x^{2n}=(2v)^{2}+(v\pm1)^2.\\
&\phantom{as}\notag
\end{align} 
}
It is clear that $(x,v)=1$ and it is not difficult to show 
that $x$ must be odd and $v$ must be even. Looking at~\eqref{eq:2} in 
$\mathbb{Z}[\,i\,]$, where $i = \sqrt{-1}$,  it can further be shown that
$v\pm1+2iv$ has to be an $n$th power and that there exists
integers $A$ and $B_{1}$ such that 
{\allowdisplaybreaks
\begin{align}\label{eq:3}
&\phantom{as}\notag\\
&v\pm1+2iv = (A+B_{1}i)^{n} \Longrightarrow  x^{2} = A^{2}+B_{1}^{2}.\\
&\phantom{as}\notag
\end{align}
}
  $A$ has to be odd and $B_{1}$ has to be even ($=2B$, say).
Suppose $n=2m+1$. Comparing real and imaginary parts of the 
first equation in \eqref{eq:3},
it follows that 
{\allowdisplaybreaks
\begin{align}\label{eq:4}
&\phantom{as}\notag\\
\pm1=B^{n}\sum_{j=0}^{m}
(-1)^j2^{2j}\left(\binom{n}{2j}\left(\frac{A}{B}\right)^{n-2j}-
\binom{n}{2j+1}\left(\frac{A}{B}\right)^{n-2j-1}\right).
\end{align}
}
Let 
{\allowdisplaybreaks
\begin{align}\label{fn}
f_{n}(x)=
\sum_{j=0}^{m}(-1)^j2^{2j}\left(\binom{n}{2j}x^{n-2j}-
\binom{n}{2j+1}x^{n-2j-1}\right).
\end{align}
}
Assume $f_{n}(x)$ is irreducible (which can easily be proved for the
cases examined: $ n = 5, 7, 11, 13$ and $17$) 
and let the roots be denoted $\theta_{1},
\cdots, \theta_{n}$. Let $\theta_{i}$ denote anyone of these roots.
  From \eqref{eq:4} and \eqref{fn} we have that
\begin{equation*}
\pm1 =\prod_{k=1}^{n}\left(A-\theta_{k}B\right) =
N_{K/Q}\left(A-\theta_{i}B\right).
\end{equation*}
Hence  $A-\theta_{i}B$  is a  unit  in  $\mathbb{Q}(\theta_{i})$.  
For the  cases  examined
 it turns out that $f_{n}(x)$ has all its roots real and that
  the rank of the group of units is
$n-1$. 
\footnote{It may be easy to show that $f_{n}(x)$ is irreducible for 
all primes $n$ and, using Sturm's Theorem, that $f_{n}(x)$
has all real roots for every prime $n$, in which case the rank of the
unit group is $n-1$ for all prime $n$. However, this is not examined here.}

Denote a set of fundamental units in $\mathbb{Q}(\theta_{i})$ by 
$\epsilon_{1}^{(i)}, \cdots, \epsilon_{n-1}^{(i)}$ and let
$\beta_{i}:=A-\theta_{i}B$. Then there exists integers 
$u_{1},\cdots, u_{n-1}$ such that
{\allowdisplaybreaks
\begin{align}\label{eq:6}
\beta_{i} = \pm \prod_{k=1}^{n-1}\epsilon_{k}^{(i)u_{k}}.
\end{align}
}
Let $U = \max_{1\leq k \leq n-1}|u_{k}|$. 
We next find an initial bound on $U$.

\section{Finding an initial bound on the exponents of the fundamental units}

Suppose $j$ is such that 
$|\beta_{j}|= \min_{1\leq i \leq n}|\beta_{i}|$ and let $K= \mathbb{Q}(\theta_{j})$. Define the following numbers:
{\allowdisplaybreaks
\begin{align}\label{ceqs}
c_{1} &= \min_{i,\,\,i \ne j}|\theta_{j} - \theta_{i}|,
&c_{2} = \max_{r \ne s \ne t \ne r}
\left|\frac{\theta_{t}-\theta_{r}}{\theta_{t}-\theta_{s}}\right|,\\ 
&\phantom{as} &\phantom{as} \notag\\
c_{2a} &= \max_{(s,t),\,\,\,j \ne s \ne t \ne j}
\left|\frac{\theta_{j}-\theta_{s}}{\theta_{j}-\theta_{t}}\right|,  
&c_{3} = \max_{i,\,\,i \ne j}|\theta_{i} - \theta_{j}|.
\phantom{asdd}\notag\\
&\phantom{as} &\phantom{as} \notag
\end{align} 
}
It is also assumed that 
{\allowdisplaybreaks
\begin{align}\label{beq1}
|B| \geq \max \left\{4, \sqrt[\uproot{2}n]{2c_{2}}\frac{2}{c_{1}}\right\}.\\
&\phantom{as}\notag
\end{align} 
}
For $i \ne j$, we have that
{\allowdisplaybreaks
\begin{align}\label{bet1}
|B||\theta_{i} - \theta_{j}|
&\leq |\beta_{j}|+|\beta_{i}|
\leq 2|\beta_{i}|
\Longrightarrow
\displaystyle{\frac{ |B|c_{1}}{2} \leq|\beta_{i} |}\\
&\phantom{as}\notag
\end{align}
}
This implies that
{\allowdisplaybreaks
\begin{align}\label{bet2}
|\beta_{j} | &= \displaystyle{\frac{1}{\prod_{i=1,i \ne j}^{n}|\beta_{i} |}
\leq \left(\frac{2}{|B|c_{1}}\right)^{n-1}}.
\end{align}
}
Similarly, it follows from \eqref{beq1} that
\begin{align}\label{bet3}
&\phantom{as}\notag\\
|\beta_{i}|\leq |B||\theta_{i} - \theta_{j}| + |\beta_{j} |
 \leq \displaystyle{
|B|\left(c_{3}+ \frac{c_{1}}{4c_{2}}\right)}=c_{4}|B|,
\\&\phantom{as}\notag
\end{align}
where\,\,$c_{4}=\displaystyle{c_{3}+c_{1}/(4c_{2})}$. 
Using  \eqref{bet1}, \eqref{beq1} and \eqref{bet3}, it follows that 
\begin{align}\label{bet4}
&\phantom{as}\notag\\
 |\log|\beta_{i} || \leq c_{5}\log|B|,
\\&\phantom{as}\notag
\end{align}
where 
$c_{5} = \displaystyle{ \max \left\{1+|\log( c_{1}/2)|/ \log 4,
1+ |\log c_{4}| /\log 4\,\,
\right\}}$.

Let $\{ i_{k}: k=1, \cdots, n-1\}$ denote the set 
$\{1,2,\cdots, n\} \setminus \{j\}$. From \eqref{eq:6},
{\allowdisplaybreaks
\begin{align*}
\left(
\begin{matrix}
\log|\beta_{i_{1}}| \\
\vdots \\
\log|\beta_{i_{n-1}}|
\end{matrix}
\right)
&=
\left(
\begin{matrix}
\log|\epsilon_{1}^{(i_{1})}| & \dots & \log|\epsilon_{n-1}^{(i_{1})}| \\
\vdots & \ddots & \vdots \\
\log|\epsilon_{1}^{(i_{n-1})}| & \dots & \log|\epsilon_{n-1}^{(i_{n-1})}|
\end{matrix}
\right)
\left(
\begin{matrix}
u_{1}\\
\vdots\\
u_{n-1}
\end{matrix}
\right)\\
&\phantom{as}\\
&:= M\left(
\begin{matrix}
u_{1}\\
\vdots\\
u_{n-1}
\end{matrix}
\right), \\
&\phantom{as}\\
\Longrightarrow \left(
\begin{matrix}
u_{1}\\
\vdots\\
u_{n-1}
\end{matrix}
\right)
&=M^{-1}\left(
\begin{matrix}
\log|\beta_{i_{1}}| \\
\vdots \\
\log|\beta_{i_{n-1}}|
\end{matrix}
\right).\\
&\phantom{as}
\end{align*}
}
Suppose  $M^{-1}:=\left(m_{r,s}\right),\,\,
1\leq r\leq n-1,\,\,\, 1\leq s\leq n-1$. Then 
{\allowdisplaybreaks
\begin{align*}
&\phantom{as}\\
U 
&\leq \max_{1 \leq r \leq n-1}\{|m_{r,1}|+\cdots+|m_{r,n-1}|\} \times
\max_{1 \leq t \leq n-1}\{|\log|\beta_{i_{t}}||\} \leq c_{6}\log|B|,\\
&\phantom{as}
\end{align*}
}
where 
$c_{6}=
\max_{1 \leq r \leq n-1}\{|m_{r,1}|+\cdots+|m_{r,n-1}|\}\times c_{5}$.
Thus 
{\allowdisplaybreaks
\begin{equation}\label{eq:8}
 \exp{(U/c_{6})} \leq |B|.
\end{equation}
}
Let
{\allowdisplaybreaks
\begin{align}\label{lam}
 \Lambda 
&=
\log\left|\frac{\theta_{j}-\theta_{k}}{\theta_{j}-\theta_{l}}\right|+
\sum_{r=1}^{n-1}u_{r}
\log\left|\frac{\epsilon_{r}^{(l)}}{\epsilon_{r}^{(k)}}\right|
=\log\left|\frac{\theta_{j}-\theta_{k}}{\theta_{j}-\theta_{l}}
\frac{A-\theta_{l}B}{A-\theta_{k}B}\right|.
\end{align}
}
 By Siegel's identity,
{\allowdisplaybreaks
\begin{align}\label{sieg}
&\phantom{as}\notag \\
\left|\frac{\theta_{l}-\theta_{k}}{\theta_{l}-\theta_{j}}
\frac{A-\theta_{j}B}{A-\theta_{k}B}\right|&=
\left|\frac{\theta_{j}-\theta_{k}}{\theta_{j}-\theta_{l}}
\frac{A-\theta_{l}B}{A-\theta_{k}B} - 1\right|.\\
&\phantom{as}\notag
\end{align}
}
From the definition of the $\beta_{i}$'s,  \eqref{ceqs},
\eqref{bet1} and \eqref{bet2}, it follows that
{\allowdisplaybreaks
\begin{align}\label{ineq1}
&\phantom{as}\notag\\
&\left|\frac{\theta_{l}-\theta_{k}}{\theta_{l}-\theta_{j}}
\frac{A-\theta_{j}B}{A-\theta_{k}B}\right| 
\leq c_{2}\left(\frac{2}{|B|c_{1}}\right)^{n}<\frac{1}{2}\\
&\phantom{as}\notag
\end{align}
}
From \eqref{sieg} and\eqref{ineq1}, it follows that
{\allowdisplaybreaks
\begin{align}\label{ineq2}
&\frac{\theta_{j}-\theta_{k}}{\theta_{j}-\theta_{l}}
\frac{A-\theta_{l}B}{A-\theta_{k}B} >\frac{1}{2}>0\\
&\phantom{as}\notag\\
&\Longrightarrow \Lambda = 
\log \frac{\theta_{j}-\theta_{k}}{\theta_{j}-\theta_{l}}
\frac{A-\theta_{l}B}{A-\theta_{k}B}\\
&\phantom{as}\notag\\
&\Longrightarrow |\Lambda|<2c_{2}\left(\frac{2}{|B|c_{1}}\right)^{n}.\\
&\phantom{as}\notag
\end{align}
}
The last inequality follows, in the case $\Lambda>0$,
 from \eqref{ineq1},  \eqref{sieg} and 
the fact that 
$e^{x}-1 > x$, for  $x>0$. In the case $\Lambda < 0$, we 
also use \eqref{ineq2}. Note that $\Lambda \not = 0$, or
else the right side of \eqref{sieg} is zero, implying,
on the left side of \eqref{sieg}, that either $\theta_{j}$ is rational,
or $\theta_{l}=\theta_{k}$. However, both of these are impossible,
since $\theta_{j}$ is algebraic of degree $n$, and 
$f_{n}(x)$ has distinct roots.

Combining this last inequality for $|\Lambda|$ with~\eqref{eq:8} it follows that 
{\allowdisplaybreaks
\begin{equation}\label{eq:10}
\frac{c_{1}^n}{2^{n+1}c_{2}} \exp{(nU/c_{6})}
<\frac{1}{|\Lambda|}.
\end{equation}
}
Next, the following theorem of Baker and W\H{u}stholz (\cite{BW93}) gives an upper bound on 
$1/|\Lambda|$\,\,:
\begin{theorem}
Denote by $\alpha\sb{1},\cdots,\alpha\sb{n}$ algebraic numbers, 
not $0$ or $1$, by
$\log\alpha\sb{1}$, $\cdots$, $\log\alpha\sb{n}$ determinations of their 
logarithms, by $d$ the degree over ${Q}$ of the number
field ${Q}(\alpha\sb{1},\cdots,\alpha\sb{n})$ and by 
$b\sb{1},\cdots,b\sb{n}$ rational integers, not all $0$
and let $B=\max\{\vert b\sb{1}\vert ,\cdots,\vert b\sb{n}\vert ,e\sp{1/d}\}$.

Define
$\log A\sb{i}= \max\{h(\alpha\sb{i}),(1/d)\vert
\log\alpha\sb{i}\vert ,1/d\}$ ($1\le i\le n$), where $h(\alpha)$ 
denotes the absolute logarithmic Weil height of $\alpha$.
   Assuming
 the number $\Lambda=b\sb{1}\log\alpha\sb{1}+\cdots+b\sb{n}\log\alpha\sb{n}$ does not vanish, then 
{\allowdisplaybreaks
\begin{align*}
\vert \Lambda\vert
\ge\exp\{-C(n,d)\log A\sb{1}\cdots\log A\sb{n}\log B\},
\end{align*}
}
 where $C(n,d)=18(n+1)!n\sp{n+1}(32d)\sp{n+2}\log(2nd)$. 
\end{theorem}
Here 
{\allowdisplaybreaks
\begin{align*}
h(\alpha)=\frac{1}{[\mathbb{Q}(\alpha):\mathbb{Q}]}
\log \left|\, a_{0}\,\prod_{r=1}^{s}\max\{1,|\alpha^{i}|\}\right|,
\end{align*}
}
where
the minimal polynomial of $\alpha$ has leading coefficient $a_{0}$ and
$\alpha = \alpha^{1},\cdots,\alpha^{s}$ are the conjugates of $\alpha$.

In our application, it can be seen from \eqref{lam} that
 $n$ has the same meaning as previously, that
 $b_{1} = 1$, that $\alpha_{1}= 
\left|(\theta_{j}-\theta_{k})/(\theta_{j}-\theta_{l})\right|$,
and, for $j=2,\cdots,n$, that $b_{j} = u_{j-1}$ and   
$\alpha_{j}=\left|\epsilon_{j-1}^{(l)}/ \epsilon_{j-1}^{(k)}\right|$.

Let $\gamma_{i} = \epsilon_{i}^{(l)}/ \epsilon_{i}^{(k)}$, with
conjugates $\gamma_{i}=\gamma_{i}^{1},\cdots,\gamma_{i}^{s}$. Since
$\gamma_{i}$ is a unit  its minimal polynomial has its leading
coefficient $a_{0} = 1$ and since $f_{n}(x)$ has all real roots,
$|\gamma_{i}^{r}|= \pm \gamma_{i}^{r}$.
{\allowdisplaybreaks
\begin{align*}
\displaystyle{|\gamma_{i}^{r}| \leq  \eta_{i} := 
\frac{\max\{|\epsilon_{i}^{(r)}|: 1 \leq r \leq n\}}
{\min\{|\epsilon_{i}^{(r)}|: 1 \leq r \leq n\}}}
 \Longrightarrow h(\gamma_{i}) \leq \log \eta_{i}.
\end{align*}
}
Let
$\delta_{jkl}
=\displaystyle{(\theta_{j}-\theta_{k})/(\theta_{j}-\theta_{l})}$ 
and suppose the minimum polynomial of $\delta_{jkl}$ is $g(x)$. 
The conjugates of $\delta_{jkl}$ are bounded by $c_{2a}$. 
For small primes $n$ the Galois group associated to the polynomial 
$f_{n}(x)$ can be determined using a computer Algebra system like
Magma.
Let $p(x)=\prod(x-\delta_{rst})(\theta_{r}-\theta_{t})$ where the product
is taken over all conjugates $\delta_{rst}$ of $\delta_{jkl}$.
$p(x)\in \mathbb{Z}[\,x\,]$ and $g(x)\,|\,p(x)$. For low values of $n$,
$p(x)$  can be calculated  numerically and  $g(x)$ and thus
$d_{g} =$ degree$(g(x))$ and $a_{0}$ can be determined explicitly.
In fact, for the cases examined, $d_{g}= n(n - 1)$ and 
$a_{0}=2^{2(n-1)}$.
Then 
{\allowdisplaybreaks
\begin{align*}
h(\delta_{jkl}) \leq \frac{\log 2^{2(n-1)}}{ n(n - 1)} + \log c_{2a}.
\end{align*}
}
{\allowdisplaybreaks
\begin{align*}
 \left[\mathbb{Q}\left(
\displaystyle{\delta_{jkl},
\epsilon_{1}^{(l)}/ \epsilon_{1}^{(k)}, \cdots,
\epsilon_{i}^{(l)}/ \epsilon_{i}^{(k)} }\right):\mathbb{Q}\right] &\leq
 [\mathbb{Q}(\theta_{1},\cdots,\theta_{n}):\mathbb{Q}] =:D.
\\&\phantom{as}
\end{align*}
}
For the values of $n$ examined, $D = n(n-1)$.

Let $C(n,D)$ be as defined in the theorem. Then
{\allowdisplaybreaks
\begin{align*}
|\Lambda|&>\exp\left(-C(n,D)\prod_{i}\log \eta_{i}
\left(\frac{\log 2^{2(n-1)}}{ n(n - 1)} + \log c_{2a}\right)
\log U\right)\\
&=\exp(-c_{7}\log U) = U^{-c_{7}},\notag
\end{align*}
}
where $c_{7}=
C(n,D)\prod_{i}\log \eta_{i}  (\frac{\log 2^{2(n-1)}}{ n(n - 1)} + \log c_{2a})$.
Combining this inequality with \eqref{eq:10} it follows that 
\begin{equation*}
\frac{c_{1}^n}{2^{n+1}c_{2}} \exp{(nU/c_{6})}<U^{c_{7}}.
\end{equation*}
If it is assumed that $|U| \geq 4$ then 
\[\frac{U}{\log U} \leq \frac{c_{6}c_{7}}{n}+ 
\frac{c_{6}\left|\log \displaystyle{\frac{c_{1}^n}{2^{n+1}c_{2}}}\right|}
{n\log 4 } .
\] 
Thus an upper bound can be found for $U$.
Denote
this upper bound by $K_{3}$. This bound is generally too large to enable the
remaining cases to be tested so it is next reduced, using the LLL
algorithm.

\section{Reducing the bound }
To reduce the bound a version of the LLL algorithm is applied, as
outlined in the paper of Tzanakis and De Weger~\cite{TD89}.
Using the notation of their paper:

Let 
\begin{equation}\label{eq:14a}
\Lambda = \delta + a_{1}\mu_{1}+ \cdots + a_{q}\mu_{q},
\end{equation} where
the $a_{i}$'s are integers and the $\mu$'s and $\delta$ are real, 
with $\delta \ne 0$. Let $A = \max_{1 \leq i \leq q}|a_{i}|$.
$K_{1},K_{2}$ and $K_{3}$ are positive numbers satisfying 
\begin{equation}\label{eq:14}
|\Lambda|<K_{1}\exp (-K_{2}A), A<K_{3}.
\end{equation}
Choose $c_{0}= \sigma_{1}K_{3}^{q}$ where $\sigma_{1}>1$.
Consider the lattice $\Gamma$  associated with the matrix 
\begin{equation*}
\mathcal{A}=
\left(
\begin{matrix}
1 &0& \dots & 0 &0\\
0&1&\dots&0&0\\
\vdots &\vdots& \ddots & \vdots & \vdots \\
0 & 0 & \dots & 1 & 0 \\
\text{[}c_{0}\mu_{1}] &[c_{0}\mu_{2}] & \dots & [c_{0}\mu_{q-1}] &  [c_{0}\mu_{q}]
\end{matrix}
\right)
\end{equation*}
Find a reduced basis $\textbf{$b_{1}$},\dots,\textbf{$b_{q}$}$ for this
basis and let $\mathcal{B}$ be the matrix associated with this basis. 
\[
\text{Let }\textbf{x}=\left(
\begin{matrix}
0 \\
\vdots \\
0\\
\text{[}-c_{0}\delta]
\end{matrix}
\right) \in \mathbb{Z}^{q}\text{ and let }\textbf{s}= \left(
\begin{matrix}
s_{1} \\
\vdots \\
s_{q-1}\\
s_{q}
\end{matrix}
\right) = \mathcal{B}^{-1}\textbf{x}.
\] Let $||y||$ denote the distance
from $y$ to the nearest integer. 
\begin{proposition}
Let $i^{*} = max\{i:1 \leq i \leq q$ and  $s_{i} \notin \mathbb{Z}\}$. If
\begin{equation*}
2^{-(q-1)/2}||s_{i^{*}}|||\textbf{$b_{1}$}|
\geq\sqrt{\left(4q^{2}+3q-\frac{3}{4}\right)}K_{3},
\end{equation*}
Then every solution of~\eqref{eq:14a} satisfying~\eqref{eq:14}
satisfies
\begin{equation*}
A<\frac{1}{K_{2}}\log \left( \frac{c_{0} K_{1}}{q K_{3}}\right). 
\end{equation*}
\end{proposition}
From \eqref{eq:10}, this proposition can be applied with $a_{r} = u_{r}$,
for  $1 \leq r \leq n-1$, and 
\begin{align*}
&q=n-1,& &A=U,& &K_{1}= \frac{2^{n+1}c_{2}}{c_{1}^n},& \\
&K_{2}= \frac{n}{c_{6}},& 
&\delta = 
\log\left|\frac{\theta_{j}-\theta_{k}}{\theta_{j}-\theta_{l}}\right|,&
&\mu_{r}=
\log\left|\frac{\epsilon_{r}^{(l)}}{\epsilon_{r}^{(k)}}\right|.& 
\end{align*} 
Once a new lower bound
is found the proposition is then applied again by now setting 
$K_{3}=1/K_{2}\log \left( c_{0} K_{1}/(q K_{3})\right)$
 and this is repeated
until the bound is reduced as far as possible. 

\section{Completing the Search}

Once $K_{3}$ has been reduced as much as the LLL algorithm will allow,
 a computer
search is done of products of powers  of the fundamental units,
with these powers bounded by
this final value of $K_{3}$,  as on
the right hand side of~\eqref{eq:6}, to see if any of these products 
have the form of the left side of~\eqref{eq:6}. However, see later
for the cases $n=11$,$13$ and $17$.

Remark: These calculations have to be carried out for each value of
$j$ -- in other words all of the $c_{i}$'s, apart from $c_{2}$ are
dependent on the choice of $j$.

These theoretical results are now applied for $n=5,7,11,13$ and $17$.

\section{Small prime powers in the Fibonacci sequence}
Remark: As noted above, most of the $c_{i}$ depend on the choice of
$j$ and hence are given as vectors (the first component being the
value got by letting $j = 1$ and so on). 
For a fixed $j$, $k$ is chosen to be $j+1 \mod{n}$ and $l$ is chosen
to be $j+2 \mod{n}$. In what follows, $\theta$ denotes a root of 
$f_{n}(x)$.

1)The case $n=5$: 
$D=20$  and 
\[f_{5}(x)=-16 + 80x + 40x^2 - 40x^3 - 5x^4 + x^5. 
\]
The zeroes of $f_{5}(x)$ are 
\[
\{ -4.64105,-1.1869,0.185992,1.75785,8.88411\}.
\]
A set of
fundamental units in $\mathbb{Q}(\theta)$ is 
{\allowdisplaybreaks
\begin{align*}
 &\{1/16\theta^3 + 3/16\theta^2 - 1/4\theta - 1/4, \\
 & 3/32\theta^4 - 27/32\theta^3 - 1/4\theta^2 + 25/8\theta - 1/2,\\ 
 & 5/16\theta^4 - 9/8\theta^3 - 223/16\theta^2 - 15/2\theta + 33/4, \\
 & 1/32\theta^4 - 13/32\theta^3 - 1/8\theta^2 + 87/8\theta - 2\}.
\end{align*}
}
{\allowdisplaybreaks
\begin{align*}
 & c_{1} = \{3.4541, 1.3728, 1.3728, 1.5718, 7.1262 \},\\
 &c_{2} = 7.3356,\\
 &c_{2a} = \{3.9156, 7.3356, 6.3356, 4.5336, 1.8979 \},\\
 &c_{3} = \{ 13.5251, 10.0710, 8.6981, 7.1262, 13.5251 \},\\
 &c_{5} = \{2.8850, 2.6694, 2.5642, 2.4219, 2.8916 \},\\
 &c_{6} = \{1.8086, 1.6734, 1.4252, 1.3461, 1.4617 \},\\
 &c_{7} = \{1.7353\times10^{32}, 2.3987\times10^{32},
 2.2438\times10^{32}, 1.89018\times10^{32}, \\
&\phantom{asfdsfsdafsadfasfsdfsadf}    9.70057\times10^{31} \},\\
 &K_{1} =  \{0.954799, 96.2577, 96.2577, 48.9276, 0.0255452\}\text{ and } \\
 &K_{2} = \{2.76452, 2.98782, 3.50828, 3.71432, 3.42053\}.
\end{align*}
}
Initially,  $K_{3} = \{10^{34},10^{34},10^{34},10^{34},10^{34}\}$,
  and eventually 
$K_{3}$  $=$  $\{11,12,$  $10,10,8\}$. 
Finally, a check on products of powers of 
fundamental units, with the power
being bounded in absolute value by $12$ shows that there are no fifth
powers in the Fibonacci sequence other than the trivial ones.

\vspace{20pt}

2)The case $n=7$: \, $D=42$ and  
\[f_{7}(x) =64 - 448 x - 336 x^2 + 560 x^3 + 140 x^4 
- 84 x^5 - 7 x^6 + x^7.
\]
The zeroes of $f_{7}(x)$ are 
\[
\{ -6.68663,-2.19286,-0.804777,0.132665,1.13197,2.88015,
  12.5395\}.
\]
 A set of fundamental units for 
$\mathbb{Q}(\theta)$ is
{\allowdisplaybreaks
\begin{align*}
\{ &5/512\theta^6 - 33/512\theta^5 - 211/256\theta^4 +
59/64\theta^3 + 31/8\theta^2 - 73/32\theta + 5/16,\\
 &  1/512\theta^6 - 7/512\theta^5 - 11/64\theta^4 + 25/64\theta^3
+ 27/32\theta^2 - 39/32\theta + 1/4, \\
 & 3/512\theta^6 - 37/512\theta^5 - 13/64\theta^4 + 177/64\theta^3
- 171/32\theta^2 + 87/32\theta +    1/4,\\
& 5/256\theta^6 - 7/64\theta^5 - 433/256\theta^4 - 9/32\theta^3 + 
69/32\theta^2 + 9/8\theta -     1/16,\\
& 7/256\theta^5 - 33/256\theta^4 - 81/32\theta^3 - 75/32\theta^2 + 
87/16\theta + 11/16, \\
 & 1/256\theta^6 - 7/256\theta^5 - 5/16\theta^4 + 11/16\theta^3 + 
39/16\theta^2 - 21/16\theta -   2\}.
\end{align*}
}
{\allowdisplaybreaks
\begin{align*}
&c_{1}=\{4.49377, 1.38808, 0.937441, 0.937441, 0.999309, 1.74818,
9.65932\},\\
&c_{2}= 14.2348,\\
&c_{2a}=\{4.27839, 10.6135, 14.2348, 13.2348, 11.4154, 5.52537,
1.99042\},\\
&c_{3}=\{19.2261, 14.7323, 13.3443, 12.4068, 11.4075, 9.65932,
19.2261\},\\
&c_{5}=\{3.13545, 2.94165, 2.86996, 2.81749, 2.75706, 2.63825,
3.13883\},\\
&c_{6}=\{2.05127, 2.13505, 2.08302, 2.04265, 1.8492, 1.61462,
2.10145\},\\
&c_{7}=\{1.03176 \times 10^{47}, 1.59932 \times 10^{47}, 
1.78271 \times 10^{47},
    1.7372 \times 10^{47},\\ &1.6448 \times 10^{47},
1.19154 \times 10^{39}, 
    5.53721 \times 10^{39}\}\\
&K_{1} = \{0.0984719, 367.017, 5727.71, 5727.71, 3661.77, 
73.0283, 0.000464476\},\\
&K_{2}= \{3.41253, 3.27862, 3.36051, 3.42691, 3.78542, 4.33539, 3.33104\}
\end{align*}
}
$K_{3}$ is initially
 $\{10^{49},10^{49},10^{49},10^{49},
10^{49},10^{49},10^{49}\}$ and  eventually
$\{16, 17,$  $ 17, 18,
 17, 15, 16\}$. A check on  products of powers of 
fundamental units, with the power
being bounded in absolute value by $18$ produces no such products 
which are linear in $\theta$ and thus 
 shows that there are no seventh
powers in the Fibonacci sequence other than the trivial ones .

Remark: A set of fundamental units for the cases $n=11,13$ and $17$ are 
not included here because they are so large. 
However, they can be found in Appendix I. They can also be easily 
generated in GP/PARI. For $n=11$, the following  code produces 
a set of fundamental units:\footnote{ Guillaume Hanrot pointed out to me that my use of the 
GP/PARI \emph{bnfinit}
command in a preprint version of this paper (where I used it without
the flag ``1'') did not necessarily produce a system of fundamental units.
Thus it is necessary to a little more to show, 
for each $n$ in question, that the system
of units I used, $\{\epsilon_{i}\}_{i=1}^{n-1}$, is indeed a 
fundamental system. Let the system of fundamental units generated by  
\emph{bnfinit}($f_{n},1$) in GP/PARI, version 2.1.0, be denoted by 
$\{\alpha_{i}\}_{i=1}^{n-1}$. 
For the cases $n=5$ and $n=7$, this  set equals 
$\{\epsilon_{i}\}_{i=1}^{n-1}$. 
For $n=11$, 
\[
\{\epsilon_{i} \}_{i=1}^{10} =
\left \{
\alpha_{1},\alpha_{2},\alpha_{3},\alpha_{4},\alpha_{5},\alpha_{8},
\alpha_{9},\alpha_{7},\frac{\alpha_{2}\alpha_{6}\alpha_{7}}
{\alpha_{3}\alpha_{4}},\alpha_{10} \right \}.
\]
For $n=13$,
\[
\{\epsilon_{i} \}_{i=1}^{12} =
\left \{
\alpha_{8},\alpha_{5},\alpha_{2},\alpha_{1},\alpha_{3},\alpha_{4},
-\frac{\alpha_{3}\alpha_{4}}
{\alpha_{6}},
\alpha_{7},
\frac{1}
{\alpha_{10}},
-\alpha_{4}\alpha_{9},
\frac{\alpha_{1}\alpha_{5}\alpha_{6}\alpha_{11}}
{\alpha_{3}\alpha_{4}\alpha_{7}},
\alpha_{12} \right \}.
\]
For $n=17$,
\begin{multline*}
\{\epsilon_{i} \}_{i=1}^{16} =
\biggm \{
\frac{\alpha_{1}^{2}\alpha_{2}\alpha_{3}^{2}\alpha_{6}\alpha_{7}^{2}
\alpha_{11}\alpha_{16}}
{\alpha_{4}^{2}\alpha_{5}\alpha_{8}^{2}\alpha_{9}\alpha_{12}},
\alpha_{2},\alpha_{4},\alpha_{3},\alpha_{5},\alpha_{6},\alpha_{7},
\alpha_{8},\alpha_{9},\alpha_{12},\\
\frac{\alpha_{1}\alpha_{2}\alpha_{3}\alpha_{11}}
{\alpha_{4}},
\alpha_{13},\alpha_{10},\alpha_{14},
-\frac{\alpha_{2}}
{\alpha_{10}\alpha_{15}},
-\frac{1}
{\alpha_{1}} \biggm \}.
\end{multline*}
Each of these sets is  clearly also a set of fundamental units.
}

\begin{align*}
\{f&= Pol(1024 - 11264x - 14080x^2 + 42240x^3 + 21120x^4\\& - 29568x^5 
- 7392x^6 +   5280x^7 + 660x^8 - 220x^9 - 11x^{10} + x^{11},x);\\
&bnfinit(f,1)[8][5]
\}
\end{align*}

\vspace{20pt}

\vspace{20pt}

3)The case $n=11$: \, $D=110$ and  
\begin{multline*}
f_{11}(x) =1024 - 11264\,x - 14080\,x^2 + 42240\,x^3 + 21120\,x^4 - 
  29568\,x^5\\ - 7392\,x^6 + 5280\,x^7 + 660\,x^8 - 
  220\,x^9 - 11\,x^{10} + x^{11}.
\end{multline*} 
The zeroes of $f_{11}(x)$ are 
\begin{multline*}
\{ -10.6902,\,-3.93193,\,-2.12056,\,-1.16928,\,-0.496752,\,0.0843495,\\
  0.680024,\,1.40783,\,2.51488,\,4.91791,\,19.8037\}.
\end{multline*}

\allowdisplaybreaks{
\begin{multline*}
c_{1}=\{6.75826,\, 1.81137,\, 0.951281,\, 0.672528,\, 0.581101,\, 0.581101,\\
 0.595674, \,
0.727806,\,1.10705,\, 2.40303,\, 14.8858\}.
\end{multline*}
$c_{2}=34.9345$.
\begin{multline*}
c_{2a}=\{4.5121,\,13.1037,\, 23.0471,\, 31.1853,\, 34.9345,\, 33.9345,\, 32.1043,\\ 25.2758,\,
 15.6171,\, 6.49517,\, 2.04852 \}.
\end{multline*}
\begin{multline*}
c_{3}=\{30.4939,\, 23.7356,\, 21.9243,\, 20.973,\, 20.3005,\, 19.7194,\, 19.1237,\\ 18.3959,\,
17.2888,\, 15.6081,\, 30.4939  \}.
\end{multline*}
\begin{multline*}
c_{5}=\{3.46637,\, 3.28489,\, 3.22745,\, 3.1954,\, 3.17187,\, 3.15092,\, 3.12881,\\ 3.10086,\,
3.05622,\, 2.98291,\, 3.46774  \}.
\end{multline*}
\begin{multline*}
c_{6}=\{3.03767,\, 2.38601,\, 3.34699,\, 3.39081,\, 2.34158,\, 3.44318,\, 3.15266,\\ 3.38847,\,
2.68879,\, 2.40831,\, 3.06079 \}.
\end{multline*}
\begin{multline*}
c_{7}=\{2.8731\times{10}^{78},\,4.7491\times{10}^{78},\,5.7427\times{10}^{78},\,
  6.2748\times{10}^{78},\\
6.4746\times{10}^{78},\, 6.4235\times{10}^{78},\,
  6.326\times{10}^{78},\,5.9051\times{10}^{78},\,\\5.0579\times{10}^{78},\,
  3.5141\times{10}^{78},\, 
1.4836\times{10}^{78}
\}
\end{multline*}
\begin{multline*}
K_{1} = \{0.0001,\,207.753,\,247867.,\,1.1241\times{10}^7,\,5.6085\times{10}^7,\,
  5.6085\times{10}^7, \\
4.271\times{10}^7,\,
4.7146\times{10}^6,\,46750.2,\,9.2739,\,1.7994\times{10}^{-8}š\}.
\end{multline*}
\begin{multline*}
K_{2}= \{3.6212,\, 4.61021,\, 3.28653,\, 3.24406,\, 4.69768,\, 3.19472,\, 3.48911,\,\\
 3.2463,\,
4.09106,\, 4.56752,\, 3.59385\}.
\end{multline*}
}
Initially,
\begin{multline*}
K_{3}= 
\{ \,{10}^{81},\,\,\,{10}^{81},\,{10}^{81},\,{10}^{81},
  \,{10}^{81},\,{10}^{81},\,\,\,{10}^{81},\,{10}^{81},
  \,{10}^{81},\,{10}^{81},\,\,    10^{80}\},
 \end{multline*}
and finally $K_{3}= \{32,\, 27,\, 41,\, 43,\, 29,\, 47,\, 40,\, 43,\, 32,\, 26,\, 28\}$. 

Checking all products of powers of ten independent units and with
these powers being bounded absolutely by $47$ 
 would take some time but there is a much quicker way, which is now
described. The same method is applied for $n=13$ and $17$.

Let $p = \sum_{i=0}^{10}a_{i}\theta^{i}$,
$q = \sum_{i=0}^{10}b_{i}\theta^{i}$ be two numbers in 
$\mathbb{Q}(\theta)$. Then $p \times q =
\sum_{i=0}^{10}c_{i}\theta^{i}$,
 for some $c_{i} \in \mathbb{Q}$ and if the $a_{i}$ are bounded in
absolute value by $K_{a}$ and the $b_{i}$ are bounded in
absolute value by $K_{b}$ then a bound for the $c_{i}$ can 
be found in terms of $K_{a}$ and $K_{b}$. In fact the 
$c_{i}$ are bounded by $16564181057933828K_{a}K_{b}$. Another way to
see this is to regard $p$ and $q$ as polynomials in $\theta$ and
reduce their product modulo $f_{11}(\theta)$. Let 
$M = 16564181057933828$.

Next, pick out the coefficient $v$ that is largest in absolute value
in the following list of powers of the fundamental units:
 $\{\epsilon_{r}^{i}: 1 \leq r \leq 10, -47 \leq i \leq 47
\}$, regarding the members of this set as polynomials in $\theta$  
(This involves checking $950$ units rather than $95^{10}$). 
We have 
{\allowdisplaybreaks
\begin{multline*}
v=1/512 \times (20107468130152762104958655475357868066478593\\
506162563506545
02987105724151326017006680926209502\\
499326050998267664485654586806568806547).
\end{multline*}
}
Thus the coefficient of any power of $\theta$ in any expression of the form
$\prod_{r=1}^{10}\epsilon_{r}^{i_{r}}$ where $-47 \leq i_{r} \leq 47$
is less than $M^{9}v^{10}$. In particular, 
$\beta_{i} = A + B\theta_{i}$ must have $A$ and $B$ less than 
$M^{9}v^{10}$ and thus, from~\eqref{eq:3},  $x \leq \sqrt{5} M^{9}v^{10}$
and thus $F_{m} = x^{11} \leq (\sqrt{5} M^{9}v^{10})^{11} 
\backsimeq 7.7943\times 10^{15864}$. Since
$F_{m} \backsimeq ((1+\sqrt{5})/2)^{m}/\sqrt{5} 
\Longrightarrow m \leq 75913$. The odd terms in the Fibonacci 
sequence up to $F_{75913}$ can be checked directly  to see if any are 
eleventh powers and no non-trivial eleventh powers are found.

More efficiently, one can use the standard trick that
 if an eleventh power exists, it must be an eleventh
power residue modulo every prime. Choose, say, ten primes
$p_{1}, \cdots, p_{10} \equiv 1
(\mod{11})$ and calculate the eleventh power residues in each case. Using the
fact that $F_{i+3} = 3F_{i+1} - F_{i-1}$ (it being necessary to
consider  $F_{j}$ for $j$ odd) and working modulo each of the ten
primes in parallel, it is a matter of seconds to check up $j=75913$
for eleventh powers (by checking if $F_{j} \mod{p_{i}}$ is an
eleventh power residue for $p_{i}$, for each $i$).

\vspace{20pt}

4)The case $n=13$: \, $D=156$ and  
\begin{multline*}
f_{13}(x) =-4096 + 53248\,x + 79872\,x^2 - 292864\,x^3\\ - 183040\,x^4 + 
  329472\,x^5 + 109824\,x^6 - 109824\,x^7\\ - 20592\,x^8 + 
  11440\,x^9 + 1144\,x^{10} - 312\,x^{11} - 13\,x^{12} + 
  x^{13}.
\end{multline*}
The roots of $f_{13}(x)$ are 
\begin{multline*}
\{ -12.6754,\,-4.75486,\,-2.68723,\,-1.64838,\,-0.960337,\,-0.41792,\\
  0.0713607,\,0.569323,\,1.14244,\,1.90337,\,3.13069,\,5.90001,\,
  23.4269\}.
\end{multline*}
\begin{multline*}
c_{1}=\{7.92054,\, 2.06763,\, 1.03885,\, 0.688043,\, 0.542417,\, 0.48928,\\ 0.48928,\, 
0.497963,\, 0.573112,\, 0.760936,\, 1.22732,\, 2.76932,\,17.52696209\,\}.
\end{multline*}
$c_{2}= 48.7346$.
\begin{multline*}
c_{2a}=\{4.55807,\, 13.63,\, 25.1375,\, 36.4444,\, 44.9604,\, 48.7346,\, 47.7346,\\ 45.9023,\,
 38.8833,\, 28.2856,\, 16.537,\, 6.70758,\, 2.05982 17.52\}.
\end{multline*}
\begin{multline*}
c_{3}=\{36.1023,\, 28.1818,\, 26.1142,\, 25.0753,\, 24.3873,\, 23.8449,\, 23.3556,\\ 22.8576,\, 
22.2845,\, 21.5236,\, 20.2963,\, 18.5754,\, 36.1023.52\}.
\end{multline*}
\begin{multline*}
c_{5}=\{3.58782,\, 3.40862,\, 3.35353,\, 3.3242,\, 3.30411,\, 3.28788,\, 3.27293,\\ 3.25738,\,
 3.23908,\, 3.21405,\, 3.17179,\, 3.10821,\, 3.588823.52\}.
\end{multline*}
\begin{multline*}
c_{6}=\{5.2067,\, 3.88525,\, 4.86669,\, 3.30557,\, 3.62425,\, 3.62461,\, 3.28336,\\ 3.98813,\,
 2.97576,\, 3.41762,\, 3.86076,\, 3.09418,\, 4.097533.52\}.
\end{multline*}
\begin{multline*}
c_{7}=\{1.30674\times{10}^{95},\,2.18839\times{10}^{95},\,2.68104\times{10}^{95},\,
  2.97999\times{10}^{95},\, \\3.14901\times{10}^{95},\,3.21389\times{10}^{95},\,
  3.1972\times{10}^{95},\,3.1657\times{10}^{95},\\
3.03213\times{10}^{95},\,
   2.77601\times{10}^{95},\,2.34399\times{10}^{95},\,\\1.6177\times{10}^{95},\,
  6.67449\times{10}^{94}\}.
\end{multline*}
\begin{multline*}
K_{1} = \{1.65365\times{10}^{-6},\,63.2569,\,486473.,\,1.03099\times{10}^8,\,
2.26948\times{10}^9, \\
8.66973\times{10}^9,\,8.66973\times{10}^9,\,
  6.89763\times{10}^9,\,1.10954\times{10}^9,\\
2.78439\times{10}^7,\, 
55693.3,\, 1.41715,\,
5.421\times{10}^{-11}\}.
\end{multline*}
\begin{multline*}
K_{2}= \{2.49678,\, 3.34599,\, 2.67122,\, 3.93275,\, 3.58695,\, 3.5866,\, 3.95935,\\ 3.25967,\, 
 4.36864,\, 3.80382,\, 3.36721,\, 4.20144,\, 3.17265
\}.
\end{multline*}
Initially,
\begin{multline*}
K_{3} 
=\{{10}^{98},\,  {10}^{98},\,  {10}^{98},\,  {10}^{98},\, 
     {10}^{98},\,  {10}^{98},\,  {10}^{98},\,  {10}^{98},\, \\
     {10}^{98},\,  {10}^{98},\,  {10}^{98},\,  {10}^{97},\, 
     {10}^{97}\}
\end{multline*}
and eventually
$K_{3}=\{53,\, 45,\, 58,\, 49,\, 55,\, 55,\, 43,\, 52,\, 44,\, 50,\, 55,\, 41,\, 47\}$.

With the same notation that was used for the case $n=11$, 
{\allowdisplaybreaks
\begin{multline*}
v=1/4096 \times (93158647867090656840416856127516852294230637148702\\
851086532124807140957259454209260273172314431029910278429059765\\
08393206322152405473550058771766196947352038793187460444181),
\end{multline*}
}
$M=316357820342343521286$, and so $F_{m} \leq (\sqrt{5}M^{11}v^{12})^{13}
\backsimeq
 5.4892\times 10^{29199}$ and so 
$m \leq 139720$. A check shows that $F_{j}$ is not a thirteenth power 
for $3 \leq j \leq 139720$.

\vspace{20pt}

5)The case $n=17$: \,$D=272$  and 
\begin{multline*}
f_{17}(x) =-65536 + 1114112\,x + 2228224\,x^2 - 11141120\,x^3 - 
  9748480\,x^4\\ + 25346048\,x^5 + 12673024\,x^6 - 
  19914752\,x^7 - 6223360\,x^8\\
   + 6223360\,x^9 + 
  1244672\,x^{10} - 792064\,x^{11} - 99008\,x^{12} +\\ 
  38080\,x^{13} + 2720\,x^{14} - 544\,x^{15} - 17\,x^{16} + 
  x^{17}.
\end{multline*}
The zeroes of $f_{17}(x)$ are
\begin{multline*}
\{ -16.6323,\,-6.36449,\,-3.75619,\,-2.50343,\,-1.72576,\,-1.16412,\\
  -0.712712,\,-0.317684,\,0.0545603,\,0.430621,\,0.838223,\,1.31512,\\
  1.92569,\,2.80429,\,4.30707,\,7.83505,\,30.6661\}.
\end{multline*}
{\allowdisplaybreaks
\begin{multline*}
c_{1}=\{10.2678,\,2.60829,\,1.25276,\,0.777669,\,0.561638,\,0.451412,\\
  0.395027,\,0.372245,\,
0.372245,\,0.376061,\,0.407602,\,0.476899,\,
  0.61057,\,\\  0.878601,\,1.50278,\, 3.52799,\,22.831
\},
\end{multline*}
$c_{2}= 83.2349$.
\begin{multline*}
c_{2a}=\{4.60646,\,14.1973,\,27.4771,\,42.6525,\,57.6738,\,70.5125,\,79.4345,\\
  83.2349,\, 82.2349,\,80.4005,\,73.1789,\,61.5455,\,47.0714,\,31.7115,\\
  17.5402,\,6.93523,\,2.07167
\}.
\end{multline*}
\begin{multline*}
c_{3}=\{47.2984,\,37.0306,\,34.4223,\,33.1695,\,32.3918,\,31.8302,\,31.3788,\\
  30.9838,\, 30.6115,\,30.2355,\,29.8279,\,29.351,\,28.7404,\,27.8618,\\
  26.359,\,24.4674,\,47.298467
\}.
\end{multline*}
\begin{multline*}
c_{5}=\{3.78233,\,3.60548,\,3.55271,\,3.52594,\,3.50882,\,3.49619,\,3.48589,\\
  3.47675,\,3.46803,\,3.45911,\,3.44932,\,3.4377,\,3.42255,\,3.40018,\\
  3.36024,\,3.30671,\,3.782967
\}.
\end{multline*}
\begin{multline*}
c_{6}=\{6.96297,\,6.95734,\,5.89133,\,6.24564,\,4.71335,\,4.94999,\,5.6139,\\
  7.93478,\,7.87795,\,4.98754,\,6.01398,\,7.84567,\,6.87894,\,5.26067,\\
  5.08113,\,5.13209,\,8.364577
\}.
\end{multline*}
\begin{multline*}
c_{7}=\{2.15293\times{10}^{126},\,3.65902\times{10}^{126},\,
  4.54254\times{10}^{126},\,5.13093\times{10}^{126},\\
  5.53464\times{10}^{126},\,5.80357\times{10}^{126},\,
  5.96299\times{10}^{126},\,6.02552\times{10}^{126},\\
  6.00935\times{10}^{126},\,5.97916\times{10}^{126},\,
  5.85324\times{10}^{126},\,5.62158\times{10}^{126},\\
  5.26283\times{10}^{126},\,
4.73432\times{10}^{126},\,
  3.94195\times{10}^{126},\,2.7004\times{10}^{126},\\
  1.08369\times{10}^{126}
\}.
\end{multline*}
\begin{multline*}
K_{1} = \{1.3922\times{10}^{-10},\,1.82303,\,473234.,\,1.56794\times{10}^9,\,
  3.9635\times{10}^{11},\\
1.62592\times{10}^{13},\,1.57104\times{10}^{14},\,
  4.3128\times{10}^{14},\,4.3128\times{10}^{14},\,3.62626\times{10}^{14},\\
  9.22206\times{10}^{13},\,6.39148\times{10}^{12},\,9.57958\times{10}^{10},\,
  1.96963\times{10}^8,\,21460.7,\\
0.01073,\,1.75352\times{10}^{-16}
\}.
\end{multline*}
\begin{multline*}
K_{2}= \{2.44149,\,2.44346,\,2.88559,\,2.7219,\,3.60678,\,3.43435,\,3.0282,\,\\
  2.14247,\,2.15792,\,
3.40849,\,2.82675,\,2.1668,\,2.47131,\,3.23152,\,
  3.34572,\,\\3.31249,\,2.03238
\}.
\end{multline*}
}
Initially, 
\begin{multline*}
K_{3}=
\{10^{134},\,\,10^{134},\,10^{134},\,10^{134},\,10^{134},\,10^{134},
10^{134},10^{134},
10^{134},\,10^{134},\,\\
10^{134},\,10^{134},\,10^{134},\,10^{134},\,10^{134},\,10^{134},\,
10^{133} \}
\end{multline*}
 and  eventually
\begin{multline*}
K_{3}=\{
  93,\,103,\, 91,\, 100,\, 76,\, 81,\, 93,\, 135,\, 134,\, 82,\, 102,\, \\
132,\, 113,\, 83,\, 77,\, 73,\, 106
  \}.
\end{multline*}
With the same notation as above $M=416654165624561667592653373446$ and 
\allowdisplaybreaks{
\begin{multline*}
v=1/16384\times (3001857560961454376246370500531976342749510086525\\3
7165233797532993409765024146699104795999297502467819521882253175\\87
914591257285898763774317774494225802733903917259383140801108235\\422
8188011341225855205786331925792324752396626091769176991951013\\90211
864167053047355724708811474915491143119130389255089524641379\\747951
73794783022459203950545394990641318064587830393038066236424\\9440973
7496287203469936364043185785741160332733990397733872532456\\24091789
0154747),
\end{multline*}
}
$F_{m} \leq (\sqrt{5}M^{15}v^{16})^{17}$ $\backsimeq 3.2504 \times
10^{128942}$ and thus 
$m \leq 616986$. A check shows that $F_{j}$ is not a seventeenth power 
in the range $3 \leq j \leq 616986$.

\section{Conclusion}
The same method could be used to extend the results about prime powers
in other binary recurrence sequences, in particular the Lucas
sequence.
It is possible the alternative method
 used to overcome the problems of have a large
number of independent units and a large bound on their exponents may
be applied in other situations. 

\vspace{10pt}

\emph{Acknowledgements:}
I wish to thank Guillaume Hanrot for several helpful comments.
He drew my attention to his paper \cite{H00}, in which he describes 
a method for solving  Thue equations without the full unit group,
a method can be applied to the problem of finding small prime
powers in the Fibonacci sequence. He also points out that,
once the LLL reduction of the bound on the powers of the units
has been completed, that there several methods for shortening the final 
check of remaining possible cases (See \cite{BH96}). 

Finally, he pointed out that the ``thueinit'' and ``thue'' commands
in PARI/GP can  be used to solve the associated
Thue equation, when the degree is small.

\section{Appendix I : Fundamental Units in $\mathbb{Q}(\theta)$ for the 
cases $n=11,13$ 
and $17$.}
For the case $n=11$ a set of fundamental units in the associated field is 
the following:
{\allowdisplaybreaks
\begin{multline*}
\biggm\{
 \frac{209}{256} - \frac{2945\,x}{256} + 
   \frac{83\,x^2}{8} + \frac{7449\,x^3}{512} - 
   \frac{21179\,x^4}{2048} - \frac{4375\,x^5}{1024} + 
   \frac{4257\,x^6}{2048} + \frac{2711\,x^7}{8192} - \\
   \frac{5987\,x^8}{65536} - \frac{323\,x^9}{65536} + 
   \frac{7\,x^{10}}{16384},\,\,\,\,
  \frac{421}{256} + \frac{1643\,x}{512} - 
   \frac{1343\,x^2}{512} - \frac{3605\,x^3}{512} - \\
   \frac{5311\,x^4}{2048}
 + \frac{4537\,x^5}{4096} + 
   \frac{2571\,x^6}{4096} - \frac{77\,x^7}{8192} - 
   \frac{1631\,x^8}{65536} - \frac{101\,x^9}{131072} + 
   \frac{13\,x^{10}}{131072},\\
  - \frac{3}{128}  + \frac{207\,x}{256} - 
   \frac{271\,x^2}{256} - \frac{39\,x^3}{256} + 
   \frac{801\,x^4}{1024} - \frac{715\,x^5}{2048} - 
   \frac{27\,x^6}{2048} + \frac{113\,x^7}{4096} - \\
   \frac{71\,x^8}{32768}
 - \frac{17\,x^9}{65536} + 
   \frac{x^{10}}{65536},\,\,\,\,- \frac{81}{256}   + 
   \frac{1147\,x}{256} + \frac{1309\,x^2}{256} - 
   \frac{2417\,x^3}{256} - \frac{5511\,x^4}{2048} + \\
   \frac{4209\,x^5}{2048} + \frac{1089\,x^6}{4096} - 
   \frac{369\,x^7}{4096} - \frac{297\,x^8}{65536} + 
   \frac{27\,x^9}{65536},\,\,\,\,
  \frac{147}{128} + \frac{47\,x}{128} + 
   \frac{1111\,x^2}{256} -\\ \frac{2577\,x^3}{256} - 
   \frac{803\,x^4}{256} + \frac{5\,x^5}{2} + 
   \frac{1353\,x^6}{4096} - \frac{463\,x^7}{4096} - 
   \frac{187\,x^8}{32768} + \frac{17\,x^9}{32768},\\
  - \frac{31}{4}   - \frac{2423\,x}{256} + 
   \frac{2271\,x^2}{128} + \frac{315\,x^3}{256} - 
   \frac{21377\,x^4}{1024} - \frac{4845\,x^5}{2048} + 
   \frac{16773\,x^6}{4096} + \frac{1491\,x^7}{4096}\\ - 
   \frac{2779\,x^8}{16384} - \frac{503\,x^9}{65536} + 
   \frac{49\,x^{10}}{65536},\,\,\,\,
  \frac{17}{16} - \frac{909\,x}{256} - 
   \frac{23\,x^2}{128} + \frac{2757\,x^3}{512} - 
   \frac{737\,x^4}{1024}\\ - \frac{1705\,x^5}{1024} - 
   \frac{55\,x^6}{4096} + \frac{659\,x^7}{8192} - 
   \frac{7\,x^8}{16384} - \frac{31\,x^9}{65536} + 
   \frac{x^{10}}{65536},\,\,\,\,
-\frac{501}{256}   + 
   \frac{2155\,x}{512} +\\ \frac{13771\,x^2}{512} + 
   \frac{2787\,x^3}{256} - \frac{49915\,x^4}{2048} - 
   \frac{23013\,x^5}{4096} + \frac{9541\,x^6}{2048} + 
   \frac{1145\,x^7}{2048} - \frac{12909\,x^8}{65536} - \\
   \frac{1273\,x^9}{131072} + \frac{117\,x^{10}}{131072},\,\,\,\,
  \frac{113}{256} + \frac{4569\,x}{256} + 
   \frac{2559\,x^2}{64} - \frac{927\,x^3}{512} - 
   \frac{61735\,x^4}{2048} - \frac{1721\,x^5}{512} + \\
   \frac{5429\,x^6}{1024} + \frac{3943\,x^7}{8192} - 
   \frac{14059\,x^8}{65536} - \frac{641\,x^9}{65536} + 
   \frac{31\,x^{10}}{32768},\\
  \frac{625}{256} - \frac{7401\,x}{256} + 
   \frac{1651\,x^2}{256} + \frac{3747\,x^3}{128} - 
   \frac{19629\,x^4}{2048} - \frac{14289\,x^5}{2048} + \\
   \frac{8857\,x^6}{4096} + \frac{439\,x^7}{1024} - 
   \frac{6255\,x^8}{65536} - \frac{365\,x^9}{65536} + 
   \frac{15\,x^{10}}{32768}
\biggm\}.
\end{multline*}
}

For $n=13$ a set of fundamental units in the associated field is:

{\allowdisplaybreaks
\normalsize{
\begin{multline*}
\biggm\{\frac{23}{1024} - \frac{10109\theta}{2048} - 
   \frac{44003\theta^2}{4096} + \frac{124839\theta^3}{8192} + 
   \frac{136209\theta^4}{4096} - \frac{20289\theta^5}{16384} - 
   \frac{433515\theta^6}{32768}\\ - \frac{61481\,\theta^7}{65536} + 
   \frac{369007\,\theta^8}{262144} + \frac{49647\,\theta^9}{524288} - 
   \frac{39319\,\theta^{10}}{1048576} - 
   \frac{2917\,\theta^{11}}{2097152} + \frac{243\,\theta^{12}}{2097152}
   ,\\ -\frac{3727}{1024}   + 
   \frac{7181\,\theta}{128} - \frac{216211\,\theta^2}{4096} - 
   \frac{121847\,\theta^3}{1024} + \frac{655867\,\theta^4}{8192} + 
   \frac{119945\,\theta^5}{2048} - \frac{991747\,\theta^6}{32768}\\ - 
   \frac{37475\,\theta^7}{4096} + \frac{888013\,\theta^8}{262144} + 
   \frac{7157\,\theta^9}{16384} - 
   \frac{101367\,\theta^{10}}{1048576} - 
   \frac{1149\,\theta^{11}}{262144} + \frac{335\,\theta^{12}}{1048576},\\
  \frac{1407}{1024} - \frac{2625\,\theta}{256} - 
   \frac{160297\,\theta^2}{4096} + \frac{2013\,\theta^3}{2048} + 
   \frac{395753\,\theta^4}{8192} + \frac{6531\,\theta^5}{1024} - 
   \frac{524081\,\theta^6}{32768}\\ - \frac{31605\,\theta^7}{16384} + 
   \frac{427491\,\theta^8}{262144} + \frac{8689\,\theta^9}{65536} - 
   \frac{45573\,\theta^{10}}{1048576} - 
   \frac{891\,\theta^{11}}{524288} + \frac{143\,\theta^{12}}{1048576},\\
  \frac{195}{256} - \frac{15321\,\theta}{2048} - 
   \frac{58735\,\theta^2}{2048} + \frac{5855\,\theta^3}{8192} + 
   \frac{289521\,\theta^4}{8192} + \frac{76579\,\theta^5}{16384} - 
   \frac{191731\,\theta^6}{16384}\\ - \frac{92625\,\theta^7}{65536} + 
   \frac{19531\,\theta^8}{16384} + \frac{50851\,\theta^9}{524288} - 
   \frac{16651\,\theta^{10}}{524288} - 
   \frac{2605\,\theta^{11}}{2097152} + \frac{209\,\theta^{12}}{2097152}
   ,\\
\frac{1}{16} - \frac{639\,\theta}{512} - 
   \frac{175\,\theta^2}{64} + \frac{3547\,\theta^3}{2048} + 
   \frac{8821\,\theta^4}{2048} - \frac{2439\,\theta^5}{4096} - 
   \frac{243\,\theta^6}{128}\\ - \frac{1129\,\theta^7}{16384} + 
   \frac{3441\,\theta^8}{16384} + \frac{1581\,\theta^9}{131072} - 
   \frac{23\,\theta^{10}}{4096} - \frac{105\,\theta^{11}}{524288} + 
   \frac{9\,\theta^{12}}{524288},\\
  \frac{243}{1024} - \frac{2135\,\theta}{1024} - 
   \frac{22321\,\theta^2}{4096} + \frac{64117\,\theta^3}{4096} + 
   \frac{53391\,\theta^4}{8192} - \frac{60195\,\theta^5}{8192} - 
   \frac{48321\,\theta^6}{32768}\\ + \frac{27909\,\theta^7}{32768} + 
   \frac{22815\,\theta^8}{262144} - \frac{6291\,\theta^9}{262144} - 
   \frac{1053\,\theta^{10}}{1048576} + \frac{81\,\theta^{11}}{1048576},\\
  \frac{243}{1024} - \frac{1659\,\theta}{512} - 
   \frac{26341\,\theta^2}{4096} + \frac{5907\,\theta^3}{512} + 
   \frac{149841\,\theta^4}{8192} - \frac{10071\,\theta^5}{4096} - 
   \frac{229029\,\theta^6}{32768}\\ - \frac{525\,\theta^7}{2048} + 
   \frac{190759\,\theta^8}{262144} + \frac{5441\,\theta^9}{131072} - 
   \frac{20001\,\theta^{10}}{1048576} - 
   \frac{89\,\theta^{11}}{131072} + \frac{61\,\theta^{12}}{1048576},\\
  - \frac{593}{512}   + \frac{19249\,\theta}{1024} - 
   \frac{51101\,\theta^2}{2048} - \frac{123517\,\theta^3}{4096} + 
   \frac{66827\,\theta^4}{2048} + \frac{124229\,\theta^5}{8192} - 
   \frac{185417\,\theta^6}{16384}\\ - \frac{85397\,\theta^7}{32768} + 
   \frac{159199\,\theta^8}{131072} + \frac{35493\,\theta^9}{262144} - 
   \frac{17713\,\theta^{10}}{524288} - 
   \frac{1529\,\theta^{11}}{1048576} + \frac{115\,\theta^{12}}{1048576}
   ,\\
\frac{797}{1024} - \frac{19295\,\theta}{2048} - 
   \frac{93697\,\theta^2}{4096} + \frac{419113\,\theta^3}{8192} + 
   \frac{307313\,\theta^4}{4096} - \frac{249571\,\theta^5}{16384} - 
   \frac{912353\,\theta^6}{32768}\\ - \frac{8383\,\theta^7}{65536} + 
   \frac{736037\,\theta^8}{262144} + \frac{67797\,\theta^9}{524288} - 
   \frac{75437\,\theta^{10}}{1048576} - 
   \frac{5115\,\theta^{11}}{2097152} + \frac{453\,\theta^{12}}{2097152}
   ,\\
- \frac{2287}{1024}   - 
   \frac{31689\,\theta}{256} - \frac{334863\,\theta^2}{4096} + 
   \frac{308411\,\theta^3}{2048} + \frac{751519\,\theta^4}{8192} - 
   \frac{11207\,\theta^5}{256} -\\ \frac{878031\,\theta^6}{32768} + 
   \frac{50093\,\theta^7}{16384} + \frac{632093\,\theta^8}{262144} + 
   \frac{813\,\theta^9}{65536} - \frac{59971\,\theta^{10}}{1048576} - 
   \frac{813\,\theta^{11}}{524288} + \frac{169\,\theta^{12}}{1048576},\\
  \frac{27}{512} + \frac{2813\,\theta}{1024} - 
   \frac{90129\,\theta^2}{2048} - \frac{427171\,\theta^3}{4096} + 
   \frac{155347\,\theta^4}{1024} + \frac{639577\,\theta^5}{8192} - 
   \frac{1076753\,\theta^6}{16384}\\ - \frac{494651\,\theta^7}{32768} + 
   \frac{983883\,\theta^8}{131072} + \frac{216065\,\theta^9}{262144} - 
   \frac{111149\,\theta^{10}}{524288} - 
   \frac{9527\,\theta^{11}}{1048576} + \frac{721\,\theta^{12}}{1048576}
   ,\\
- \frac{303}{512}   + \frac{909\,\theta}{128} + 
   \frac{4857\,\theta^2}{2048} - \frac{5137\,\theta^3}{128} + 
   \frac{71615\,\theta^4}{4096} + \frac{24109\,\theta^5}{1024} - 
   \frac{146907\,\theta^6}{16384}\\ - \frac{7933\,\theta^7}{2048} + 
   \frac{136005\,\theta^8}{131072} + \frac{5593\,\theta^9}{32768} - 
   \frac{16083\,\theta^{10}}{524288} - \frac{25\,\theta^{11}}{16384} + 
   \frac{55\,\theta^{12}}{524288}\biggm\}.
\end{multline*}
}
}
Finally, for the case $n=17$ a set of fundamental units in the associated 
field is:
{\allowdisplaybreaks
\normalsize{
\begin{multline*}
\biggm\{ \frac{1669}{8192} + \frac{1013441\,\theta}{32768} - 
   \frac{831099\,\theta^2}{16384} - 
   \frac{15553837\,\theta^3}{131072} + 
   \frac{9705537\,\theta^4}{65536} + 
   \frac{70985341\,\theta^5}{524288}\\ - 
   \frac{33065551\,\theta^6}{262144} - 
   \frac{124377673\,\theta^7}{2097152} + 
   \frac{21555899\,\theta^8}{524288} + 
   \frac{90616595\,\theta^9}{8388608} - 
   \frac{22540321\,\theta^{10}}{4194304}\\ - 
   \frac{26767543\,\theta^{11}}{33554432} + 
   \frac{4417175\,\theta^{12}}{16777216} + 
   \frac{2772335\,\theta^{13}}{134217728} - 
   \frac{256021\,\theta^{14}}{67108864} - 
   \frac{66067\,\theta^{15}}{536870912} + \\
   \frac{3807\,\theta^{16}}{536870912},\\
  - \frac{3381}{16384}   + 
   \frac{1443\,\theta}{512} - \frac{1104995\,\theta^2}{65536} - 
   \frac{596807\,\theta^3}{32768} + 
   \frac{12341635\,\theta^4}{262144} + 
   \frac{1742777\,\theta^5}{65536}\\ - 
   \frac{39959399\,\theta^6}{1048576} - 
   \frac{6905549\,\theta^7}{524288} + 
   \frac{50345477\,\theta^8}{4194304} + 
   \frac{1358191\,\theta^9}{524288} - 
   \frac{25774521\,\theta^{10}}{16777216}\\ - 
   \frac{1694321\,\theta^{11}}{8388608} + 
   \frac{4980657\,\theta^{12}}{67108864} + 
   \frac{91457\,\theta^{13}}{16777216} - 
   \frac{285757\,\theta^{14}}{268435456} - 
   \frac{4507\,\theta^{15}}{134217728} + \\
   \frac{527\,\theta^{16}}{268435456},\\
  \frac{13}{64} - \frac{11119\,\theta}{4096} - 
   \frac{5847\,\theta^2}{512} + \frac{359545\,\theta^3}{16384} + 
   \frac{899457\,\theta^4}{16384} - \frac{810279\,\theta^5}{65536} - \\
   \frac{1690541\,\theta^6}{32768} - \frac{148823\,\theta^7}{262144} + 
   \frac{4347743\,\theta^8}{262144} + 
   \frac{1161099\,\theta^9}{1048576} - 
   \frac{544193\,\theta^{10}}{262144} - \\
   \frac{616077\,\theta^{11}}{4194304} +
   \frac{406671\,\theta^{12}}{4194304} + 
   \frac{87947\,\theta^{13}}{16777216} - 
   \frac{11269\,\theta^{14}}{8388608} - 
   \frac{2581\,\theta^{15}}{67108864} + \\
   \frac{161\,\theta^{16}}{67108864},\\
  \frac{3177}{4096} - \frac{37351\,\theta}{8192} - 
   \frac{260767\,\theta^2}{16384} + \frac{854595\,\theta^3}{32768} + 
   \frac{5165539\,\theta^4}{65536} - \frac{160099\,\theta^5}{131072} - \\
   \frac{17391923\,\theta^6}{262144} - 
   \frac{5968849\,\theta^7}{524288} + 
   \frac{18808385\,\theta^8}{1048576} + 
   \frac{6236715\,\theta^9}{2097152} - 
   \frac{8469197\,\theta^{10}}{4194304} - \\
   \frac{2025815\,\theta^{11}}{8388608} + 
   \frac{1527105\,\theta^{12}}{16777216} + 
   \frac{217439\,\theta^{13}}{33554432} - 
   \frac{85281\,\theta^{14}}{67108864} - 
   \frac{5339\,\theta^{15}}{134217728} + \\
   \frac{313\,\theta^{16}}{134217728},\\
  \frac{35177}{4096} + \frac{87353\,\theta}{8192} - 
   \frac{1685455\,\theta^2}{16384} - \frac{1286573\,\theta^3}{32768} + 
   \frac{18808867\,\theta^4}{65536} + 
   \frac{10349213\,\theta^5}{131072}\\ - 
   \frac{62922435\,\theta^6}{262144} - 
   \frac{26769697\,\theta^7}{524288} + 
   \frac{79984865\,\theta^8}{1048576} + 
   \frac{25086155\,\theta^9}{2097152} - 
   \frac{40730397\,\theta^{10}}{4194304}\\ - 
   \frac{8793575\,\theta^{11}}{8388608} + 
   \frac{7786305\,\theta^{12}}{16777216} + 
   \frac{1029535\,\theta^{13}}{33554432} - 
   \frac{441361\,\theta^{14}}{67108864} - 
   \frac{26891\,\theta^{15}}{134217728} + \\
   \frac{1609\,\theta^{16}}{134217728},\\
  \frac{8353}{16384} + \frac{82469\,\theta}{32768} - 
   \frac{334519\,\theta^2}{65536} - \frac{2389817\,\theta^3}{131072} + 
   \frac{4045291\,\theta^4}{262144} + 
   \frac{16465649\,\theta^5}{524288} - \\
   \frac{21174147\,\theta^6}{1048576} - 
   \frac{41242981\,\theta^7}{2097152} + 
   \frac{41290081\,\theta^8}{4194304} + 
   \frac{33247615\,\theta^9}{8388608} - 
   \frac{25999301\,\theta^{10}}{16777216} -\\ 
   \frac{9888859\,\theta^{11}}{33554432} + 
   \frac{5556713\,\theta^{12}}{67108864} + 
   \frac{993323\,\theta^{13}}{134217728} - 
   \frac{335425\,\theta^{14}}{268435456} - 
   \frac{22599\,\theta^{15}}{536870912} + \\
   \frac{1269\,\theta^{16}}{536870912},\\
  \frac{17529}{8192} - \frac{199385\,\theta}{16384} - 
   \frac{2570675\,\theta^2}{32768} - \frac{1099275\,\theta^3}{65536} + 
   \frac{22040999\,\theta^4}{131072} + 
   \frac{14557163\,\theta^5}{262144} - \\
   \frac{67044783\,\theta^6}{524288} - 
   \frac{35241615\,\theta^7}{1048576} + 
   \frac{82863837\,\theta^8}{2097152} + 
   \frac{30445541\,\theta^9}{4194304} - 
   \frac{41990713\,\theta^{10}}{8388608} - \\
   \frac{10034929\,\theta^{11}}{16777216} + 
   \frac{8051613\,\theta^{12}}{33554432} + 
   \frac{1124425\,\theta^{13}}{67108864} - 
   \frac{458917\,\theta^{14}}{134217728} - 
   \frac{28469\,\theta^{15}}{268435456} + \\
   \frac{1683\,\theta^{16}}{268435456},\\
  -\frac{15439}{16384}   + 
   \frac{7457\,\theta}{4096} - \frac{772229\,\theta^2}{65536} - 
   \frac{225567\,\theta^3}{4096} + \frac{19970281\,\theta^4}{262144} + 
   \frac{5384407\,\theta^5}{65536} - \\
   \frac{83962553\,\theta^6}{1048576} - 
   \frac{5203139\,\theta^7}{131072} + 
   \frac{118428135\,\theta^8}{4194304} + 
   \frac{7938383\,\theta^9}{1048576} - 
   \frac{63912463\,\theta^{10}}{16777216} - \\
   \frac{149561\,\theta^{11}}{262144} + 
   \frac{12691267\,\theta^{12}}{67108864} + 
   \frac{249817\,\theta^{13}}{16777216} - 
   \frac{739187\,\theta^{14}}{268435456} - 
   \frac{2983\,\theta^{15}}{33554432} + \\
   \frac{1375\,\theta^{16}}{268435456},\\
  \frac{37669}{16384} - \frac{782207\,\theta}{32768} - 
   \frac{3931563\,\theta^2}{65536} + 
   \frac{29154715\,\theta^3}{131072} + 
   \frac{38848887\,\theta^4}{262144} - 
   \frac{120152419\,\theta^5}{524288} - \\
   \frac{105613751\,\theta^6}{1048576} + 
   \frac{158737183\,\theta^7}{2097152} + 
   \frac{108365877\,\theta^8}{4194304} - 
   \frac{77583405\,\theta^9}{8388608} - 
   \frac{44555329\,\theta^{10}}{16777216} + \\
   \frac{12880225\,\theta^{11}}{33554432} + 
   \frac{6875085\,\theta^{12}}{67108864} - 
   \frac{438065\,\theta^{13}}{134217728} - 
   \frac{313501\,\theta^{14}}{268435456} - 
   \frac{7771\,\theta^{15}}{536870912} + \\
   \frac{913\,\theta^{16}}{536870912},\\
  \frac{26149}{8192} - \frac{89167\,\theta}{32768} - 
   \frac{84531\,\theta^2}{2048} + \frac{2581855\,\theta^3}{131072} + 
   \frac{6955001\,\theta^4}{65536} - 
   \frac{13565819\,\theta^5}{524288} - \\
   \frac{9794367\,\theta^6}{131072} + 
   \frac{15629499\,\theta^7}{2097152} + 
   \frac{5367431\,\theta^8}{262144} - 
   \frac{1891597\,\theta^9}{8388608} - 
   \frac{2435655\,\theta^{10}}{1048576} - \\
   \frac{2428259\,\theta^{11}}{33554432} + 
   \frac{1707267\,\theta^{12}}{16777216} + 
   \frac{541463\,\theta^{13}}{134217728} - 
   \frac{45335\,\theta^{14}}{33554432} - 
   \frac{19063\,\theta^{15}}{536870912} + \\
   \frac{1259\,\theta^{16}}{536870912},\\
  -\frac{91977}{8192}   + 
   \frac{530479\,\theta}{16384} + \frac{729419\,\theta^2}{32768} - 
   \frac{6980111\,\theta^3}{65536} + 
   \frac{1869453\,\theta^4}{131072} + 
   \frac{27472683\,\theta^5}{262144} - \\
   \frac{18142681\,\theta^6}{524288} - 
   \frac{42686171\,\theta^7}{1048576} + 
   \frac{31187031\,\theta^8}{2097152} + 
   \frac{28182221\,\theta^9}{4194304} - 
   \frac{18441695\,\theta^{10}}{8388608} - \\
   \frac{7619597\,\theta^{11}}{16777216} + 
   \frac{3855551\,\theta^{12}}{33554432} + 
   \frac{725657\,\theta^{13}}{67108864} - 
   \frac{232019\,\theta^{14}}{134217728} - 
   \frac{15945\,\theta^{15}}{268435456} + \\
   \frac{883\,\theta^{16}}{268435456},\\
  - \frac{7709}{16384}   + 
   \frac{71959\,\theta}{8192} - \frac{663207\,\theta^2}{65536} - 
   \frac{282043\,\theta^3}{8192} + \frac{9339035\,\theta^4}{262144} + 
   \frac{5176181\,\theta^5}{131072} - \\
   \frac{34425699\,\theta^6}{1048576} - 
   \frac{4589823\,\theta^7}{262144} + 
   \frac{46980309\,\theta^8}{4194304} + 
   \frac{6714145\,\theta^9}{2097152} - 
   \frac{25183029\,\theta^{10}}{16777216} - \\
   \frac{61627\,\theta^{11}}{262144} + 
   \frac{5002729\,\theta^{12}}{67108864} + 
   \frac{202099\,\theta^{13}}{33554432} - 
   \frac{292097\,\theta^{14}}{268435456} - 
   \frac{2379\,\theta^{15}}{67108864} + \\
   \frac{545\,\theta^{16}}{268435456},\\
  \frac{310213}{16384} + \frac{379055\theta}{32768} - 
   \frac{10426491\theta^2}{65536} - 
   \frac{10145735\theta^3}{131072} + 
   \frac{89578131\theta^4}{262144} + 
   \frac{62976811\theta^5}{524288} - \\
   \frac{273648223\theta^6}{1048576} - 
   \frac{136720787\theta^7}{2097152} + 
   \frac{336238489\theta^8}{4194304} + 
   \frac{116558733\theta^9}{8388608} - 
   \frac{169272961\theta^{10}}{16777216} - \\
   \frac{38784309\,\theta^{11}}{33554432} + 
   \frac{32283889\,\theta^{12}}{67108864} + 
   \frac{4403545\,\theta^{13}}{134217728} - 
   \frac{1832805\,\theta^{14}}{268435456} - 
   \frac{112817\,\theta^{15}}{536870912} + \\
   \frac{6703\,\theta^{16}}{536870912},\\
  - \frac{3375}{8192}   + 
   \frac{59225\,\theta}{4096} + \frac{2751061\,\theta^2}{32768} + 
   \frac{3266747\,\theta^3}{32768} - 
   \frac{12989039\,\theta^4}{131072} - \frac{222733\,\theta^5}{1024} - \\
   \frac{24349275\,\theta^6}{524288} + 
   \frac{41990993\,\theta^7}{524288} + 
   \frac{71338227\,\theta^8}{2097152} - 
   \frac{9003067\,\theta^9}{1048576} - 
   \frac{43379513\,\theta^{10}}{8388608} + \\
   \frac{1560157\,\theta^{11}}{8388608} + 
   \frac{8421299\,\theta^{12}}{33554432} + 
   \frac{43097\,\theta^{13}}{8388608} - 
   \frac{454593\,\theta^{14}}{134217728} - 
   \frac{10569\,\theta^{15}}{134217728} + 
   \frac{383\,\theta^{16}}{67108864},\\
  \frac{17285}{16384} - \frac{465327\,\theta}{32768} - 
   \frac{4632383\,\theta^2}{65536} + \frac{899987\,\theta^3}{131072} + 
   \frac{48004271\,\theta^4}{262144} + 
   \frac{15671757\,\theta^5}{524288} - \\
   \frac{152102419\,\theta^6}{1048576} - 
   \frac{52147849\,\theta^7}{2097152} + 
   \frac{187354341\,\theta^8}{4194304} + 
   \frac{51864899\,\theta^9}{8388608} - 
   \frac{93440845\,\theta^{10}}{16777216} - \\
   \frac{18753527\,\theta^{11}}{33554432} + 
   \frac{17634053\,\theta^{12}}{67108864} + 
   \frac{2246239\,\theta^{13}}{134217728} - 
   \frac{992033\,\theta^{14}}{268435456} - 
   \frac{59731\,\theta^{15}}{536870912} + \\
   \frac{3601\,\theta^{16}}{536870912},\\
  \frac{3503}{8192} - \frac{76383\,\theta}{8192} - 
   \frac{770045\,\theta^2}{32768} + \frac{733515\,\theta^3}{16384} + 
   \frac{4622523\,\theta^4}{131072} - 
   \frac{5753349\,\theta^5}{131072} -\\ 
   \frac{8650317\theta^6}{524288} + \frac{472931\theta^7}{32768} + 
   \frac{6723545\theta^8}{2097152} - 
   \frac{3876569\theta^9}{2097152} - 
   \frac{2151455\theta^{10}}{8388608} + 
   \frac{369863\theta^{11}}{4194304} + \\
   \frac{246705\,\theta^{12}}{33554432} - 
   \frac{41315\,\theta^{13}}{33554432} - 
   \frac{7063\,\theta^{14}}{134217728} + 
   \frac{67\,\theta^{15}}{33554432} + \frac{3\,\theta^{16}}{134217728}\biggm\}.
\end{multline*}
}
}

\end{document}